\documentclass[12pt]{amsart}
\usepackage{amssymb,latexsym,esint}
\usepackage{enumerate}
\usepackage{float}
\usepackage{import}

\usepackage{mathtools}
\usepackage{amsfonts}
\usepackage{amssymb}
\usepackage{amsmath}
\usepackage{graphicx}
\usepackage{graphicx,color}
\usepackage[normalem]{ulem}
\usepackage{todonotes}
\usepackage{mathrsfs}
\usepackage{import}
\usepackage{xifthen}
\usepackage{pdfpages}
\usepackage{transparent}
\usepackage{yfonts}

\numberwithin{equation}{section}
\makeatletter
\@namedef{subjclassname@2020}{%
  \textup{2020} Mathematics Subject Classification}
\makeatother

\usepackage{amsmath,amstext,amsthm,amsxtra,mathtools,amssymb}

\usepackage[bookmarksopen,bookmarksdepth=3,colorlinks,citecolor=red,pagebackref,hypertexnames=true]{hyperref}
\usepackage[backrefs,msc-links,nobysame,initials,non-sorted-cites]{amsrefs}

\usepackage[normalem]{ulem}
\usepackage{lmodern}

\usepackage{lmodern}% http://ctan.org/pkg/lm

\usepackage[backrefs,msc-links,nobysame,initials]{amsrefs}

\allowdisplaybreaks

\newtheorem{thm}{Theorem}
\newtheorem{cor}{Corollary}
\newtheorem{lem}{Lemma}

\theoremstyle{definition}

\theoremstyle{remark}

\begin{document}

\title[$L^p(\mathbb{R}^2)$ bounds for geometric maximal operators ]{$L^p(\mathbb{R}^2)$ bounds for geometric maximal operators associated to homothecy invariant convex bases}

%\author{Dmitry Dmitrishin}
%\address{D. D.: Department of Applied Mathematics, Odessa National Polytechnic University, Odessa 65044, Ukraine}
%\email{\href{mailto: dmitrishin@opu.ua}{dmitrishin@opu.ua}}

\author{Paul Hagelstein}
\address{P. H.: Department of Mathematics, Baylor University, Waco, Texas 76798}
\email{\href{mailto:paul_hagelstein@baylor.edu}{paul\_hagelstein@baylor.edu}}
\thanks{P. H. is partially supported by a grant from the Simons Foundation (\#521719 to Paul Hagelstein).}

\author{Alexander Stokolos}
\address{A. S.: Department of Mathematical Sciences, Georgia Southern University, Statesboro, Georgia 30460}
\email{\href{mailto:astokolos@GeorgiaSouthern.edu}{astokolos@GeorgiaSouthern.edu}}


%%%%%%%%%%%%%%%%%%%%%%%%%%%%%% SECTION SECTION SECTION
\begin{abstract}
Let $\mathcal{B}$ be a nonempty homothecy invariant collection of convex sets  of positive finite  measure in $\mathbb{R}^2$.   Let $M_\mathcal{B}$ be the geometric maximal operator defined by
$$M_\mathcal{B}f(x) = \sup_{x \in R \in \mathcal{B}}\frac{1}{|R|}\int_R |f|\;.$$  We show that either $M_\mathcal{B}$ is bounded on $L^p(\mathbb{R}^2)$ for every $1 < p \leq \infty$ or that $M_\mathcal{B}$ is unbounded on $L^p(\mathbb{R}^2)$ for every $1 \leq p < \infty$.
As a corollary, we have that any density basis that is a homothecy invariant collection of convex sets in $\mathbb{R}^2$  must differentiate $L^p(\mathbb{R}^2)$ for every $1 < p \leq \infty$.
\end{abstract}

\maketitle
 
\section{Introduction}

Let $\mathcal{B}$ denote a nonempty collection of sets of positive measure in $\mathbb{R}^n$.   Associated to $\mathcal{B}$ is the \emph{geometric maximal operator} $M_\mathcal{B}$ defined on measurable functions $f$ on $\mathbb{R}^n$ by
$$M_\mathcal{B}f(x) = \sup_{x \in R \in \mathcal{B}}\frac{1}{|R|}\int_R |f|\;,$$  
where the supremum is over all members $R$ in $\mathcal{B}$ containing $x$.
If $\mathcal{B}$ is the collection of all balls in $\mathbb{R}^n$, then $M_{\mathcal{B}}$ is the well-known \emph{Hardy-Littlewood maximal operator} $M_{HL}$.   If $\mathcal{B}$ is the collection of all rectangular parallelepipeds in $\mathbb{R}^n$ whose sides are parallel to the coordinate axes, then $M_\mathcal{B}$ is the \emph{strong maximal operator} $M_S$.   If $\mathcal{B}$ is the collection of all rectangles in $\mathbb{R}^2$ whose longest sides have slope of the form $2^{-k}$ for some natural number $k$, then $M_\mathcal{B}$ is the \emph{lacunary maximal operator} $M_{lac}$.  If $\mathcal{B}$ is the collection of all rectangular parallelepipeds in $\mathbb{R}^n$, then $M_{\mathcal{B}}$ is the \emph{Kakeya-Nikodym maximal operator} which we denote here by $M_{KN}$.
\\

The $L^p$ boundedness properties of  these particular maximal operators are now well understood.  $M_{HL}$ and $M_{S}$ are bounded on $L^p(\mathbb{R}^n)$ for $1 < p \leq \infty$, $M_{lac}$ is bounded on $L^p(\mathbb{R}^2)$ for $1 < p \leq \infty$, but $M_{KN}$ is unbounded on $L^p(\mathbb{R}^n)$ for all $1 \leq p < \infty$.   For proofs of these results we refer the reader to \cite{stein2}.
 \\
 
 $M_{\mathcal{B}_\Omega}$ is said to be a \emph{directional maximal operator} if there is a nonempty set of directions $\Omega$ for which $\mathcal{B}_\Omega$ consists of every rectangle oriented in one of the directions in $\Omega$.  Bateman proved in \cite{bateman} that if $M_{\mathcal{B}_\Omega}$ is a directional maximal operator associated to directions in the plane, then either $M_{\mathcal{B}_\Omega}$ is bounded on $L^p(\mathbb{R}^2)$ for all $1 < p \leq \infty$ or $M_{\mathcal{B}_\Omega}$ is unbounded on $L^p(\mathbb{R}^2)$ for all $1 \leq p < \infty$.  This type of dichotomy for directional maximal operators on the plane was in many respects anticipated by earlier papers including \cite{batemankatz} and \cite{stok1995}. Such an explicit dichotomy is still unknown for directional maximal operators acting on functions on $\mathbb{R}^n$ for $n \geq 3$, although progress has been made on this problem by, among others, Parcet and Rogers \cite{parcetrogers} and Kroc and Pramanik \cite{krocpramanik}.
 \\
 
 It is natural to consider the question of, if $\mathcal{B}$ is a homothecy invariant collection of convex sets in $\mathbb{R}^n$ (not necessarily consisting of  \emph{all} rectangles or rectangular parallelepipeds oriented in a certain set of directions), whether $M_\mathcal{B}$ must be bounded on $L^p(\mathbb{R}^n)$ for all $1 < p \leq \infty$ or whether it must be unbounded on $L^p(\mathbb{R}^n)$ for all $1 \leq p < \infty$.   This is a difficult problem.   The purpose of this paper is to provide a solution in the two-dimensional case.    Our solution builds on ideas accumulated over the past four decades associated to the behavior of maximal operators associated to a lacunary set of directions, the Besicovitch and Kakeya conjectures, Bernouilli percolation, and sticky maps.  Particularly relevant to our solution are papers of Nagel, Stein, and Wainger \cite{nsw};  Sj\"ogren and Sj\"olin \cite{sjsj}; Lyons \cite{lyons}; Katz \cite{katz1996}; Katz, \L aba, and Tao \cite{ktl2000}; Bateman and Katz \cite{batemankatz}; and Bateman \cite{bateman}.

 \begin{thm}\label{t1}
 Let $\mathcal{B}$ be a nonempty homothecy invariant collection of convex sets of positive finite measure in $\mathbb{R}^2$.   Let $M_\mathcal{B}$ be the geometric maximal operator defined by
$$M_\mathcal{B}f(x) = \sup_{x \in R \in \mathcal{B}}\frac{1}{|R|}\int_R |f|\;.$$   Either  $M_\mathcal{B}$ is bounded on $L^p(\mathbb{R}^2)$ for every $1 < p \leq \infty$ or  $M_\mathcal{B}$ is unbounded on $L^p(\mathbb{R}^2)$ for every $1 \leq p < \infty$.
 \end{thm}

 An immediate application of this theorem falls in the classical topic of differentiation of integrals. 
 Let $\mathcal{B}$ be a homothecy invariant collection of bounded sets in $\mathbb{R}^2$ of positive measure.  Recall that $\mathcal{B}$ is said to be a \emph{density basis} provided, given a measurable set $E \subset \mathbb{R}^2$, for \mbox{a.e. $x$} we have
$\lim_{k \rightarrow \infty}\frac{1}{|R_k|}\int_{R_k} \chi_E = \chi_E(x)$
whenever $\{R_k\}$ is a sequence of sets in $\mathcal{B}$ \mbox{containing $x$} whose diameters are tending to 0.   Hagelstein and Stokolos proved in \cite{hs} that a density basis consisting of convex sets in $\mathbb{R}^2$ must differentiate $L^p(\mathbb{R}^2)$ for sufficiently large $p$.   As a corollary of Theorem 1, we then have the following.

\begin{cor}\label{cor1}
Let $\mathcal{B}$ be a homothecy invariant collection of convex sets of positive finite measure in $\mathbb{R}^2$.   If $\mathcal{B}$ is a density basis, then $M_\mathcal{B}$ is bounded on $L^p(\mathbb{R}^2)$ for $1 < p \leq \infty$ and $\mathcal{B}$  must differentiate $L^p(\mathbb{R}^2)$ for $1 < p \leq \infty$.
\end{cor}
 
 It is worth noting that the Hagelstein-Stokolos proof in \cite{hs} that a homothecy invariant density basis consisting of convex sets must differentiate $L^p(\mathbb{R}^n)$ for sufficiently large $p$ falls in the classical realm of differentiation of integrals, relying on delicate arguments involving covering lemmas and the Calder\'on-Zygmund decomposition.  The proof of Corollary \ref{cor1}, however, relies on not only these classical ideas in differentiation theory but also techniques that are Fourier analytic in nature.  In many respects this illustrates a paradigm in contemporary harmonic analysis that $L^p$ bounds for maximal operators such as the lacunary maximal operator should in theory be provable via covering lemmas (see A. C\'ordoba and R. Fefferman \cite{cf}), but at the moment the best known bounds (in particular in the range $1 < p \leq 2$) have been attained only via Fourier analytic arguments involving Littlewood-Paley theory, square functions, and complex interpolation.
 \\

 The remainder of this paper will be devoted to a proof of Theorem \ref{t1}.  Our proof will rely heavily on the ideas of Bateman and Katz \cite{bateman, batemankatz} related to recognizing the fact that the boundedness of a directional maximal operator $M_{\mathcal{B}_\Omega}$ on $L^p(\mathbb{R}^2)$ corresponds to the ability to cover $\Omega$ by finitely many $N$-lacunary sets; we also take advantage of recent ideas of Gauvan \cite{gauvan} that enable us to modify the Bateman-Katz methodology to accommodate scenarios when $\mathcal{B}$ is a homothecy invariant basis of convex sets not necessarily corresponding to a directional maximal operator.
 \\
 
 We wish to thank the referee for  insightful comments and suggestions regarding this paper.
 \section{Proof of Theorem \ref{t1}}
 
 \begin{proof}[Proof of Theorem \ref{t1}]
 If $S \subset \mathbb{R}^2$, and $x \in \mathbb{R}^2$, we define the translate $\tau_x S$ of $S$ to be the set in $\mathbb{R}^2$ satisfying
 $$\chi_{\tau_x S} (y) = \chi_S (y - x)\;.$$   If $r > 0 $, we define the dilate $r S$ of $S$ to be the set in $\mathbb{R}^2$ satisfying
 $$\chi_{r S} (y) = \chi_S (y / r)\;.$$
 
 Recall that if $\mathcal{B}$ is a homothecy invariant basis of sets in $\mathbb{R}^2$, then for any $R \in \mathcal{B}$, $r > 0$, and $x \in \mathbb{R}^2$ we have $\tau_x R \in \mathcal{B}$ and $rR \in \mathcal{B}$.
 \\

Now,  \emph{since the members of $\mathcal{B}$ are convex}, by a result of Lassak \cite{lassak} we have that if $S \in \mathcal{B}$ then there exists a rectangle $R$ containing $S$ such that a translate of a $\frac{1}{2}$-fold dilate of $R$ is contained in $S$.  
Hence we may assume without loss of generality that every set $S \in \mathcal{B}$ is covered by a rectangle $R_S$ so that $|R_S| \leq 4 |S|$, a translate of a $\frac{1}{2}$-fold  dilate of $R_S$ is contained in $S$, and such that a longest side of $R_S$ has slope in the set $[0,1]$.   

Given a set $E \subset \mathbb{R}^2$, let $\pi_x (E) = \left\{x \in \mathbb{R} : (x,y)\in E \textup{ for some } y \in \mathbb{R}\right\}.$   Let  $\pi_y (E) = \left\{y \in \mathbb{R} : (x,y)\in E \textup{ for some } x \in \mathbb{R}\right\}.$ 
Let $S \in \mathcal{B}$ be given.   Let $R_S$ be a rectangle associated to $S$ as above.   Associate to $R_S$ a parallelogram $P_S$ containing $R_S$ whose left and right hand sides are parallel to the $y$-axis, such that the slope of the bottom side is greater than or equal to 0 and strictly less than 1, such that the length of the left hand side is $2^{-k}$ times the length of $\pi_x(P_S)$ for some nonnegative integer $k$, such that the length of $\pi_y(P_S)$  is an integer multiple of the length of the left hand side, such that $|P_S| \leq 32 |R_S|$, and such that a translate of a $\frac{1}{32}$-fold dilate of $P_S$ is contained in $R_S$.  Figure \ref{fig:figure1} provides a depiction of a convex set $S$ with its associated rectangle $R_S$ and parallelogram $P_S$.  

\begin{figure}[ht]
 \centering % centering figure
 \scalebox{0.5} % rescale the figure by a factor of 0.8
 {\includegraphics{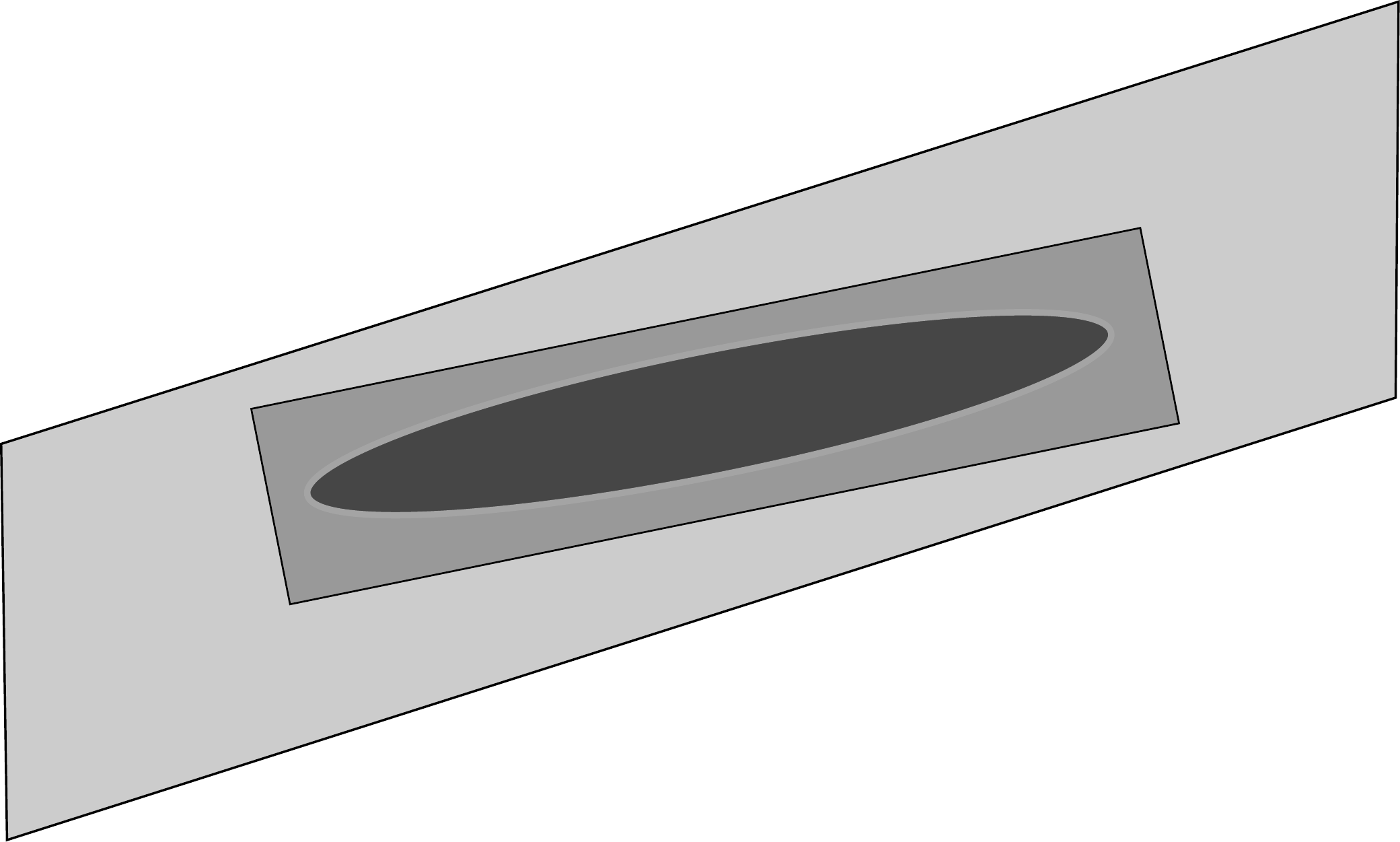}} % importing figure
 \caption{A Convex Set $S$ with Associated Rectangle $R_S$ and Parallelogram $P_S$}
\label{fig:figure1} % labeling to refer it inside the text
\end{figure}

 Let now $\tilde{\mathcal{P}}_{\mathcal{B}}$ consists of the set of all the homothecies of the parallelograms in $\{P_S : S \in \mathcal{B}\}$.   For any measurable function $f$ on $\mathbb{R}^2$ we have
$$ \frac{1}{4096} M_{\tilde{\mathcal{P}}_{\mathcal{B}}}f(x) \leq M_\mathcal{B}f(x) \leq 128 M_{\tilde{\mathcal{P}}_{\mathcal{B}}}f(x)\;.$$
Accordingly, we recognize that it suffices to show that either the maximal operator $M_{\tilde{\mathcal{P}}_{\mathcal{B}}}$ is bounded on $L^p(\mathbb{R}^2)$ for all $1 < p \leq \infty$ or that $M_{\tilde{\mathcal{P}}_{\mathcal{B}}}$ is unbounded on $L^p(\mathbb{R}^2)$ for all $1 \leq p < \infty$.
\\

 We denote by  $\mathcal{P}_{\mathcal{B}}$ the set of parallelograms in $\tilde{\mathcal{P}}_\mathcal{B}$  whose vertices are at $(0,0)$, $(0,2^{-k})$, $(1, j 2^{-k})$, and $(1, (j+1)2^{-k})$, where $k\geq 0$ is an integer and $j$ is an integer satisfying $0 \leq j \leq 2^k - 1$.  % As we are assuming without loss of generality that the slope of all the parallelograms in $\mathcal{P}_{\mathcal{B}}$ are between 0 and 1, we recognize that any parallelogram in $\mathcal{P}_{\mathcal{B}}$ is contained in a translate of a 2-fold dilate of a parallelogram in $\tilde{\mathcal{P}}_{\mathcal{B}}$ and vice versa.  Hence it suffices to show that either the maximal operator $M_{\mathcal{P}_{\mathcal{B}}}$ is bounded on $L^p(\mathbb{R}^2)$ for all $1 < p \leq \infty$ or that $M_{\mathcal{P}_{\mathcal{B}}}$ is unbounded on $L^p(\mathbb{R}^2)$ for all $1 \leq p < \infty$.
We now associate to the collection $\mathcal{P}_{\mathcal{B}}$ a subset of the dyadic tree, following the ideas of Gauvan in \cite{gauvan}.  In doing so we need to introduce relevant terminology and notation.
\\

 Let $\mathscr{B}_0 = \{0\}$, let $\mathscr{B}_1 = \{00, 01\}$, and let $\mathscr{B}_2 = \{000, 001, 010, 011\}.$   Continuing in this way, we can define the sets $\mathscr{B}_n$ recursively by $$\mathscr{B}_{k+1} = \left\{0a_1\cdots a_{k+1} : 0a_1 \cdots a_k \in \mathscr{B}_k \textup{ and } a_{k+1} \in \{0,1\}\right\}\;.$$  Let $\mathscr{B} = \cup_{n=0}^\infty \mathscr{B}_n$.  Note  $\mathscr{B}$ consists of all finite strings of 0's and 1's such that the first number in the string is 0 .  We endow $\mathscr{B}$ with the structure of a tree by adding edges between any element $0a_1 \cdots a_k$ of $\mathscr{B}$ and each of the two elements  $0a_1 \cdots a_k 0$ and $0a_1 \cdots a_k 1$.
 
 If $v \in \mathscr{B}_k$, we  say $v$ is of \emph{height} $k$ and write $h(v) = k$.

 Let $\mathscr{T}$ be a subtree of $\mathscr{B}$ .     A \emph{ray} in $\mathscr{T}$  is a (possibly infinite) maximal sequence of connected vertices $p_1, p_2, \ldots$ in $\mathscr{T}$ such that there exists a nonnegative integer $k$ so that $p_1 \in \mathscr{B}_k$, $p_2 \in \mathscr{B}_{k+1}$, $\ldots$.  It is maximal in the sense that if the ray $R$ is a finite sequence of vertices $p_1, p_2, \ldots, p_N$, then there is no vertex in $\mathscr{T}$ of height $h(p_N) + 1$ that is connected to $p_N$.  If a ray in $\mathscr{B}$ exists that contains vertices $u$ and $v$ where $h(u) < h(v)$, we say $u$ is an \emph{ancestor} of $v$ and that $v$ is a \emph{descendant} of $u$.
 
  We say that a vertex $v \in \mathscr{T}$ \emph{splits in } $\mathscr{T}$ if $v$ has two descendants in $\mathscr{T}$ of height $v(h) + 1$.
  
   We denote the set of rays in $\mathscr{T}$ whose vertex of lowest height is $v$ by $\mathfrak{R}_\mathscr{T}(v)$.   If $R$ is a ray in $\mathscr{T}$ that contains exactly $k$ vertices that split in $\mathscr{T}$, we say that $\textup{split}_\mathscr{T}(R) = k$.   If $R$ contains infinitely many vertices that split in $\mathscr{T}$, we say $\textup{split}_\mathscr{T}(R) = \infty$.
  If $v$ is a vertex in $\mathscr{T}$, we define $$\textup{split}_\mathscr{T}(v) = \min_{R \in \mathfrak{R}_\mathscr{T}(v)}\textup{split}_\mathscr{T}(R)\;.$$   We set
 $$\textup{split}(v, \mathscr{T}) = \sup_{\mathscr{S} \subset \mathscr{T}} \textup{split}_\mathscr{S}(v)\;,$$ where the supremum is over all subtrees $\mathscr{S}$ of $\mathscr{T}$,  and $$\textup{split}(\mathscr{T}) = \sup_{v \in \mathscr{T}} \textup{split}(v, \mathscr{T})\;.$$
 
 A tree $\mathscr{T} \subset \mathscr{B}$ is said to be \emph{lacunary of order 0} if $\mathscr{T}$ is a ray in $\mathscr{B}$.    A tree $\mathscr{T} \subset \mathscr{B}$ is said to be \emph{lacunary of order $N$} if all the vertices of $\mathscr{T}$ that split in $\mathscr{T}$ lie on a lacunary tree of order $N-1$.
 
We define the \emph{truncation} $\mathscr{S}^k$ of the set $\mathscr{S} \subset \mathscr{B}$ to be the  subset of $\mathscr{S}$ whose vertices consist of the vertices of $\mathscr{S}$ of height less than or equal to $k$.
 
We may associate to any  subset $\mathscr{S} \subset \mathscr{B}$  the associated  subgraph $[\mathscr{S}]$ of $\mathscr{B}$ that is the smallest subtree of $\mathscr{B}$ containing  all the vertices of $\mathscr{S}$ \emph{together with all of their ancestors in $\mathscr{B}$}.   Moreover, we may associate to $\mathscr{S}$ its \emph{extension} $\mathscr{S}^\ast$ which we define to be the induced subgraph of $\mathscr{B}$ consisting of all the vertices of $[\mathscr{S}]$ together with all vertices of the form $0a_1 a_2 a_3 \cdots a_k 0 \cdots 0$ where $0a_1 a_2 a_3 \ldots a_k \in [\mathscr{S}]$.
 \\
 
 Given the collection of parallelograms $\mathcal{P}_\mathcal{B}$ let $\mathscr{S}_\mathcal{B}$ denote the  subset of  $\mathscr{B}$ consisting of 0 and  vertices of the form
 $a_0 a_1 a_2 \ldots a_k$ (note here $a_0$ is always 0) where $\mathcal{P}_\mathcal{B}$ contains the parallelogram denoted by $P_{a_0 a_1 \cdots a_k}$ with vertices $$(0,0), (0,2^{-k}), (1, \sum_{j=0}^k a_j 2^{-j}), (1,  \sum_{j=0}^k a_j 2^{-j} + 2^{-k})\;.$$
 
\begin{lem}\label{l1}  Suppose, given $N > 0$, there exists $k > 0$ such that  $\textup{split}[\mathscr{S}_{\mathcal{B}}^k] = N$.  Then $M_\mathcal{B}$ is not bounded on $L^p(\mathbb{R}^2)$ for any $1 \leq p < \infty$.   
 \end{lem}
 
 \begin{proof}
We assume without loss of generality that  the unit square $[0,1] \times [0,1]$ lies in $\mathcal{P}_\mathcal{B}$.    Fix $N > 0$ and let $k$ be such that  $\textup{split}[\mathscr{S}_{\mathcal{B}}^k] =N$.    Assume without loss of generality that $k$ is the minimum number satisfying this condition, and in particular that \mbox{$\textup{split}[\mathscr{S}_{\mathcal{B}}^{k-1}] =N-1$.} Following the pruning technique of Bateman [p. 62 of \cite{bateman}] we prune the tree $[\mathscr{S}_\mathcal{B}^k]$ to yield a tree $\mathscr{P}_\mathcal{B}^k$ containing $v_0 = 0$ with $\textup{split}(v_0, \mathscr{P}_\mathcal{B}^k) = N$ satisfying the condition that, for every $R \in \mathfrak{R}_{\mathscr{P}_\mathcal{B}^k}(v_0)$ and every $j = 1, \ldots, N$, $R$ contains exactly one splitting \mbox{vertex $v_j$}  such that $\textup{split}(v_j, \mathscr{P}_\mathcal{B}^k) = j$. 

Let now $\tilde{\mathscr{P}}_\mathcal{B}^k$ be an induced subgraph of $\mathscr{B}$ consisting of all of the vertices in $\mathscr{P}_\mathcal{B}^k$ together with all vertices of the form $a_0 a_1 \cdots a_n 0\cdots 0$ in $\cup_{j=0}^k \mathscr{B}_j$ where $a_0a_1 \cdots a_n$ is a member of $\mathscr{P}_\mathcal{B}^k$ that has no descendant in $\mathscr{P}_\mathcal{B}^k$.   %Readers familiar with the terminology of Gauvan in \cite{gauvan} may recognize that $\tilde{\mathscr{P}}_\mathcal{B}^k$ is a fig tree of scale $N$ and height $k+1$.
\\

Let $\sigma: \mathscr{B}^k \rightarrow \tilde{\mathscr{P}}_\mathcal{B}^k$ be a \emph{sticky map}.   In particular, we have that $\sigma(u)$ is a descendant of $\sigma(v)$ in $\tilde{\mathscr{P}}_\mathcal{B}^k$ whenever $u$ is a descendant of $v$ in $\mathscr{B}^k$. We also suppose  that $h(\sigma(v)) = h(v)$ for all $v \in \mathscr{B}^k$.   Now, we can associate to the  map $\sigma$ the set $K_\sigma \subset \mathbb{R}^2$ defined by

$$K_\sigma = \bigcup_{v \in \mathscr{B}_k} E_{\sigma(v)}\;,$$
where, if $v = b_0 b_1 \cdots b_k$ and if $P_{\sigma(v)} = P_{a_0 a_1 \cdots a_k}$ is as defined earlier, then $E_{\sigma(v)}$ is the parallelogram with vertices at $$(0, \sum_{j=0}^k b_j 2^{-j}), (0,\sum_{j=0}^k b_j 2^{-j} + 2^{-k}), (2, \sum_{j=0}^k b_j 2^{-j} + 2\sum_{j=0}^k a_j 2^{-j}), (2,  \sum_{j=0}^k b_j 2^{-j} + 2\sum_{j=0}^k a_j 2^{-j} + 2^{-k})\;.$$

Now, Bateman observed that there exists such a sticky map $\sigma$ so that

$$\left|K_\sigma \cap \left([0,1] \times \mathbb{R}\right)\right| \gtrsim \frac{\log N}{N}$$
and 
$$\left|K_\sigma \cap \left([1,2] \times \mathbb{R}\right)\right| \lesssim \frac{1}{N}\;.$$

Let us abbreviate the set $K_\sigma \cap \left([0,1] \times \mathbb{R}\right)$ as $K_1$ and the set $K_\sigma \cap \left([1,2] \times \mathbb{R}\right)$ as $K_2$.
\\

It turns out that $M_{\tilde{\mathcal{P}}_\mathcal{B}} \chi_{K_2} (x)\geq \frac{1}{4}$ for every $x \in K_1$.  To see this, note that if $x \in K_1$ then $x \in \tau_{(0, j_x\cdot 2^{-k})}P_{a_0 \cdots a_k}$ for some $a_0 \cdots a_k$ in $\tilde{\mathscr{P}}_\mathcal{B}^k$, with $j_x \cdot 2^{-k} = \sum_{n=0}^{k}b_n 2^{-n}$ and $\sigma(b_0 \cdots b_k) = a_0 \cdots a_k$. %(In other words, $x$ lies in the parallelogram $P_{a_0\cdots a_k}$  translated vertically upward a distance $j_x\cdot2^{-k}$.) 
 Of course, the parallelogram $P_{a_0 \cdots a_k}$ might not lie in $\mathcal{P}_\mathcal{B}$ itself.   However, let $a_0 \cdots a_l$ be the nearest ancestor in $\mathscr{P}_\mathcal{B}^k$ to the vertex $a_0 \cdots a_k$. (We do allow $a_0 \cdots a_l$ to be $a_0 \cdots a_k$ itself provided $a_0 \cdots a_k \in \mathscr{S}_\mathcal{B}$.)  \emph{Note that} $a_0 \cdots a_k$ \emph{is the only descendant of $a_0 \cdots a_l$ in $\mathscr{P}_\mathcal{B}^k$ of height $k$.} To see this, recognize that either $a_0 \cdots a_k$ lies in the tree $\mathscr{P}_\mathcal{B}^k$ itself, or it is of the form $a_0 \cdots a_l 0 \cdots 0$ where $a_0 \cdots a_l$ is a leaf of the tree $\mathscr{P}_\mathcal{B}^k$ (i.e., a member of $\mathscr{P}_\mathcal{B}^k$ having no descendant in $\mathscr{P}_\mathcal{B}^k$), noting that all of the descendants in $\tilde{\mathscr{P}}_\mathcal{B}^k$ of a leaf $a_0 \cdots a_l$ in $\mathscr{P}_\mathcal{B}^k$ are of the form $a_0 \cdots a_l 0 \cdots 0$.   One should also recognize that if $a_0 \cdots a_l$ is a leaf of $\mathscr{P}_\mathcal{B}^k$,  then $a_0 \cdots a_l$ must lie in $\mathscr{S}_\mathcal{B}$ itself.   It is at this stage that we take advantage of the fact that $\sigma$ is a sticky map.    Since $\sigma$ is a sticky map, we must have that $\sigma(b_0 \cdots b_s) = a_0 \cdots a_s$ for $l \leq s \leq k$, and hence $\sigma(v) = a_0 \cdots a_k$ whenever $v$ is a descendant of height $k$ of $b_0 \cdots b_l$ in $\mathscr{B}_k$.   This tells us that $x$ is contained in a union of $2^{k-l}$ a.e. disjoint parallelograms whose union forms a parallelogram  $P_x$ contained in $K_1$ that is a vertical translate of $P_{a_0 \cdots a_l} \in \mathcal{P}_\mathcal{B}$ contained in $K_1$.  Not only that, but $x$ is contained in   a vertical translate $E_x$ of the associated parallelogram $E_{a_0 \cdots a_l}$ that is contained in $K_\sigma$.   As $ E_x$ is contained in a translate of a two-fold dilate of $P_x \in \tilde{\mathcal{P}}_\mathcal{B}$, we have that $M_{\tilde{P}_\mathcal{B}}\chi_{K_2}(x) \geq \frac{1}{4}$.   In particular, we have that
$$\left|\left\{x \in \mathbb{R}^2 : M_{\tilde{\mathcal{P}}_\mathcal{B}}\chi_{K_2}(x) \geq \frac{1}{4}\right\}\right| \geq |K_1| \gtrsim \log N \left|K_2\right|\;.$$  As $N > 0$ is arbitrary, the lemma holds.

 \end{proof}

\begin{lem} \label{l2} Suppose that there exists $N > 0$ so that for every $k > 0$ we have $\textup{split}[\mathscr{S}_{\mathcal{B}}^k] \leq N$.   Then $M_\mathcal{B}$ is bounded on $L^p(\mathbb{R}^2)$ for all $1 < p \leq \infty$. 
 \end{lem}
 
\begin{proof}
If for every $k > 0$  we have $\textup{split}[\mathscr{S}_{\mathcal{B}}^k] \leq N$, then $\mathscr{S}_\mathcal{B}^\ast$ is $N+1$ lacunary. 
 Let  $\Omega_{\mathcal{B}^\ast}$ denote the set of all directions associated to the $N+1$ lacunary tree  $\mathscr{S}_\mathcal{B}^\ast$.  (This could formally be defined as the smallest closed set in [0,1] containing all values of the form $\sum_{j=0}^k a_j 2^{-j}$ where $a_0 \cdots a_k \in \mathscr{S}_\mathcal{B}^\ast$.)   Denote the associated directional maximal operator by $M_{\Omega_{\mathcal{B}^\ast}}$.   As $\Omega_{\mathcal{B}^\ast}$ may be covered by a finite number of $N+1$ lacunary sets, we have by Bateman's Theorem that the directional maximal operator $M_{\Omega_{\mathcal{B}^\ast}}$   is bounded on $L^p(\mathbb{R}^2)$ for all $1 < p \leq \infty$.  As $M_{\tilde{\mathcal{P}}_\mathcal{B}} f(x) \leq  M_{\Omega_{\mathcal{B}^\ast}}f(x)$,  we have $M_\mathcal{B}$ is bounded on $L^p(\mathbb{R}^2)$ for all $1 < p \leq \infty$.
 \end{proof}
 
 As, given $\mathcal{B}$ the hypotheses of one of Lemma \ref{l1} or Lemma \ref{l2} must hold, we have proven the desired theorem.

\end{proof}

\end{document}